\documentclass[3p,times]{elsarticle}

\usepackage{ecrc}

\volume{00}

\firstpage{1}

\runauth{D. Huang et al.}

\usepackage{amssymb}

\usepackage[figuresright]{rotating}
\newcommand{\ba}{\begin{array}}
\newcommand{\ea}{\end{array}}
\newcommand{\babc}{\begin{abc}}
\newcommand{\eabc}{\end{abc}}
\newcommand{\bc}{\begin{center}}
\newcommand{\ec}{\end{center}}
\newcommand{\be}{\begin{equation}}
\newcommand{\ee}{\end{equation}}
\newcommand{\bea}{\begin{eqnarray}}
\newcommand{\eea}{\end{eqnarray}}
\newcommand{\beas}{\begin{eqnarray*}}
\newcommand{\eeas}{\end{eqnarray*}}
\newcommand{\bh}{\begin{hangitem}}
\newcommand{\eh}{\end{hangitem}}
\newcommand{\bhi}{\begin{hangitem}}
\newcommand{\ehi}{\end{hangitem}}
\newcommand{\bi}{\begin{itemize}}
\newcommand{\ei}{\end{itemize}}
\newcommand{\bn}{\begin{enumerate}}
\newcommand{\en}{\end{enumerate}}
\newcommand{\bq}{\begin{quote}}
\newcommand{\eq}{\end{quote}}
\newcommand{\btb}{\begin{tabular}}
\newcommand{\etb}{\end{tabular}}
\newcommand{\eqdef}{\stackrel{\triangle}{=}} 
\def\litem[#1]{\item[#1\hfill]}         
\newtheorem{theorem}{Theorem}

\newtheorem{lemma}{Lemma}

\newtheorem{remark}{Remark}

\def\bfa{{\bf a}}

\def\bff{{\bf f}}
\def\bfg{{\bf g}}

\def\bfu{{\bf u}}

\def\bfw{{\bf w}}
\def\bfx{{\bf x}}

\begin{document}

\begin{frontmatter}



\dochead{}

\title{Long-time Average Cost Control of Polynomial Systems: A Sum-of-squares-based Small-feedback Approach \tnoteref{t1}}


\author[im]{Deqing Huang\corref{cor1}}\ead{d.huang@imperial.ac.uk},    
\author[im]{Sergei Chernyshenko}\ead{s.chernyshenko@imperial.ac.uk}

\cortext[cor1]{Corresponding author}
\tnotetext[t1]{Funding from EPSRC under the grant EP/J011126/1 and support
in kind from Airbus Operation Ltd., ETH Zurich (Automatic Control Laboratory), University of Michigan (Department of Mathematics), and
University of California, Santa Barbara (Department of Mechanical Engineering) are gratefully acknowledged.}

\address[im]{Department of Aeronautics, Imperial College London, Prince Consort Road, London SW7 2AZ, United Kingdom}

\begin{abstract}
The two main contributions of this paper are a proof of concept of the recent novel idea in the area of long-time average cost control, and a new method of overcoming the well-known difficulty of non-convexity of simultaneous optimization of a control law and an additional tunable function.
A recently-proposed method of obtaining rigorous bounds of long-time average cost is first outlined for the uncontrolled system with polynomials of system state on the right-hand side. In this method the polynomial constraints are relaxed to be sum-of-squares and formulated as semi-definite programs.   It was proposed to use the upper bound of long-time average cost as the objective function instead of the time-average cost itself in controller design. In the present paper this suggestion is implemented for a particular system and is shown to give good results.
Designing the optimal controller by this method requires optimising simultaneously both the control law and a tunable function similar to the Lyapunov function. The new approach proposed and implemented in this paper for overcoming the inherent non-convexity of this optimisation is based on a formal assumption that the amplitude of control is small. By expanding the tunable function and the bound in the small parameter, the long-time average cost is reduced by minimizing the respective bound in each term of the series. The derivation of all the polynomial coefficients in controller is given in terms of the solvability conditions of state-dependent linear and bilinear inequalities. The resultant sum-of-squares problems are solved in sequence, thus avoiding the non-convexity in optimization. The proposed approach is implemented for a simple model of oscillatory vortex shedding behind a cylinder.
\end{abstract}

\begin{keyword}
Sum of squares; Long-time average; Polynomial systems; Small feedback; Non-convexity



\end{keyword}

\end{frontmatter}

\section{Introduction}


Although global stabilization of dynamical systems is of importance in system theory and engineering \cite{Kh:02,An:02},
it is sometimes difficult or impossible to synthesize a global stabilizing controller for certain linear and nonlinear systems \cite{Va:02}.
The reasons could be the poor controllability of system, e.g., systems that have uncontrollable linearizations \cite{Di:11}
and systems that have fewer degrees of control freedom than the degrees of freedom to be controlled \cite{Gu:13,Gu:14},
the input/output constraints in practice, e.g., an unstable linear time-invariant system cannot be globally stabilized in the presence of
input saturations \cite{Bl:99}, time delay \cite{Sun:11,Sun:13}, and/or the involved large disturbances \cite{Ki:06}, etc.
Moreover, in many applications the full stabilization, while possible, carries high penalty due to the cost of the control, thus is also not desirable.

Instead, minimizing a long-time average of the cost functional might be more realistic.
For instance, long-time-average cost analysis and control is often considered in irrigation, flood control, navigation, water supply, hydroelectric power, computer communication networks, and other applications \cite{Du:10,Bo:92}.
In addition, systems that include stochastic factors are often controlled in the sense of long-time average.
In \cite{Ro:83}, a summary of long-time-average cost problems for continuous-time Markov processes is given.
In \cite{Me:00}, the long-time-average control of a class of problems that arise in the modeling of semi-active suspension systems was considered,
where the cost includes a term based on the local time process diffusion.
Notice that the controller design methods proposed in \cite{Ro:83,Me:00} are highly dependent on the stochastic property of dynamical systems.

In certain cases, as, for example, turbulent flows of fluid, calculating the time averages is a big challenge even in the uncontrolled case.  As a result, developing the control aimed at reducing the time-averaged cost for turbulent flows, for example by using the receding horizon technique, leads to controllers too complicated for practical implementation  \cite{Bewley:01}.
To overcome this complexity, it was proposed \cite{Ph:14} to use an upper bound for the long-time average cost instead of the long-time average cost itself in cases when such an upper bound is easier to calculate. The idea is based on the hope that the control reducing an upper bound for a quantity will also reduce the quantity itself. Meanwhile, \cite{Ph:14} uses the sum of squares (SOS) decomposition of polynomials and semidefinite programming (SDP) and allows a trade-off between the quality of bound and the complexity of its calculation.

The SOS methods apply to systems defined by a polynomial vector field. Such systems may describe a wide variety of dynamics~\cite{Va:01} or approximate a system defined by an analytical vector field~\cite{Va:02}. A polynomial system can therefore yield a reliable model of a dynamical system globally or in larger regions than the linear approximation in the state-space~\cite{Va:03}. Recent results on SOS decomposition have transformed the verification of non-negativity of polynomials into SDP, hence providing promising algorithmic procedures for stability analysis of polynomial systems. However, using SOS techniques for optimal control, as for example in~\cite{Pr:02,Zh:07,Ma:10}, is subject to a generic difficulty: while the problem of optimizing the candidate Lyapunov function certifying the stability for a closed-loop system for a given controller and the problem of optimizing the controller for a given candidate Lyapunov function are reducible to an SDP and thus, are tractable, the problem of simultaneously  optimizing both the control and the Lyapuniov function is non-convex. Iterative procedures were proposed for overcoming this difficulty~\cite{Zh:07,Zh:09,Ng:11}.

While optimization of an upper bound with control proposed in~\cite{Ph:14} does not involve a Lyapunov function, it does involve a similar tunable function, and it shares the same difficulty of non-convexity. In the present work we propose a polynomial type state feedback controller design scheme for the long-time average upper-bound control, where the controller takes the structure of an asymptotic series in a small-amplitude perturbation parameter.
By fully utilizing the smallness of the perturbation parameter, the resultant SOS optimization problems are solved in sequence,
thus avoiding the non-convexity in optimization. We apply it to an illustrative example and demonstrate that it does allow to reduce the long-time average cost even without fully stabilizing the system. Notice the significant conceptual difference between our approach and the studies of control by small perturbations, often referred to as tiny feedback, see for example~\cite{tc:93}.

The paper is organized as follows. Section~\ref{seq:Background} 
  presents some preliminary introduction on SOS and its application in bound estimation of long-time average cost for uncontrolled systems. Section~\ref{seq:Problem Formulation} 
gives the problem formulation. Bound optimization of the long-time average cost for controlled polynomial systems is considered in Section~\ref{seq:Bound}.
 An illustrative example of a cylinder wake flow is addressed in Section~\ref{seq:Example}. 
Section~\ref{seq:Conclusion} 
 concludes the work.

\section{Background}\label{seq:Background}

In this section SOS of polynomials and a recently-proposed method of obtaining rigorous bounds of long-time average cost via SOS for uncontrolled polynomial systems are introduced.

\subsection{SOS of polynomials}

SOS techniques have been frequently used in the stability analysis and controller design for all kinds of systems, e.g., constrained ordinary differential equation systems \cite{An:02}, hybrid systems \cite{An:05}, time-delay systems \cite{An:04}, and partial differential equation systems \cite{Pa:06,Yu:08,GC:11}. These techniques help to overcome the common drawback of approaches based on Lyapunov functions: before \cite{Pr:02}, there were no coherent and tractable computational methods for constructing Lyapunov functions.

A multivariate polynomial $f({\bf x})$ is a SOS, if there exist polynomials $f_1({\bf x}), \cdots, f_m({\bf x})$ such that
\[
f({\bf x})=\sum_{i=1}^mf_i^2({\bf x}).
\]
If $f({\bf x})$ is a SOS then $f({\bf x})\ge 0, \forall{\bf x} $. In the general multivariate case, however,
$f({\bf x})\ge 0~  \forall \bf x$  does not necessarily imply that $f({\bf x})$ is SOS. While being stricter, the condition that $f({\bf x})$ is SOS is much more computationally tractable than non-negativity \cite{Par:00}. At the same time, practical
experience indicates that in many cases replacing non-negativity with the SOS property leads to satisfactory results.

In the present paper we will utilize the existence of efficient numerical methods and software \cite{Pra:04,Lo:09} for solving the optimization problems of the following type: minimize the linear objective function
\begin{equation}
{\bf w}^T{\bf c}
\label{linear}
\end{equation}
where $\bfw$ is the vector of weighting coefficients for the linear objective function, and
 ${\bf c}$ is a vector formed from the (unknown) coefficients of the
polynomials $p_i({\bf x})$ for $i=1,2,\cdots, \hat{N}$ and SOS $p_i({\bf x})$ for $i=(\hat{N}+1),\cdots, {N}$,
such that
\bea
&a_{0,j}({\bf x})+\sum_{i=1}^Np_i({\bf x})a_{i,j}({\bf x})=0, ~ j=1,2,\cdots, \hat{J}, \label{c5} \\
[1ex]
&a_{0,j}({\bf x})+\sum_{i=1}^Np_i({\bf x})a_{i,j}({\bf x})  \mbox{~are SOS,~} j=(\hat{J}+1),\cdots, {J}.~
\label{c6}
\eea
In (\ref{c5}) and (\ref{c6}), the $a_{i,j}(\bfx)$ are given scalar constant coefficient polynomials.

The lemma below that provides a sufficient condition to test inclusions of sets defined by polynomials is frequently used for feedback controller design
in Section~~\ref{seq:Bound}. It is a particular case of the 
Positivstellensatz 
 Theorem \cite{Po:99} and is a generalized ${\mathcal S}$-procedure \cite{Ta:06}.

\begin{lemma}
Consider two sets of $\bfx$,
\beas
{\mathcal S}_1&\eqdef& \left\{\bfx\in {\mathbb R}^n~|~h(\bfx)=0, f_1(\bfx)\ge 0, \cdots, f_r(\bfx)\ge 0\right\}, \\
[1ex]
 {\mathcal S}_2&\eqdef& \left\{\bfx\in {\mathbb R}^n~|~f_0(\bfx)\ge 0\right\},
\eeas
where $f_i(\bfx), i=0,\cdots, r$ and $h(\bfx)$ are scalar polynomial functions. The set containment ${\mathcal S}_1\subseteq {\mathcal S}_2$ holds if there exist a polynomial function $m(\bfx)$ and SOS polynomial functions $S_i(\bfx), i=1,\cdots, r$ such that
\beas
f_0(\bfx)-\sum_{i=1}^rS_i(\bfx)f_i(\bfx)+m(\bfx)h(\bfx) ~~\mbox{is SOS}.
\eeas
\end{lemma}

\subsection{Bound estimation of long-time average cost for uncontrolled systems}

For the convenience of the reader we outline here the method of obtaining bounds for long-time averages proposed in~\cite{Ph:14} and make some remarks on it.
Consider a system
\bea
\dot{\bfx}=\bff(\bfx),
\label{sys1}
\eea
where $\dot{\bfx}\eqdef d\bfx/dt$ and $\bff(\bfx)$ is a vector of multivariate polynomials of the components of the state vector $\bfx\in{\mathbb R}^n$.  The long-time average of a function of the state $\Phi(\bfx)$ is defined as
\[
\bar{\Phi}=\lim_{T\rightarrow \infty}\frac{1}{T}\int_0^T\Phi({\bfx}(t))\,dt,
\]
where $\bfx(t)$ is the solution of (\ref{sys1}).

Define a polynomial function of the system state, $V(\bfx)$, of degree $d_V,$ and containing unknown decision variables as its coefficients.
The time derivative of $V$ along the trajectories of system (\ref{sys1}) is
\[
    \dot{V}(\bfx) = \dot{\bfx} \cdot \nabla_{\bfx} V(\bfx)=\bff(\bfx)\cdot \nabla_{\bfx} V(\bfx).
 \]
Consider the following quantity:
\[
    H(\bfx) \eqdef \dot{V}(\bfx) + \Phi(\bfx)= \bff(\bfx) \cdot\nabla_{\bfx} V(\bfx) + \Phi(\bfx).
\]
The following result is from \cite{Ph:14}:

\begin{lemma}
For the system (\ref{sys1}), assume that the state $\bfx$ is bounded in $\mathcal{D}\subseteq {\mathbb R}^n$. Then, $ H(\bfx)\le C, \forall \bfx\in \mathcal{D}$ implies $\bar{\Phi}\le C$.
\label{lemma1}
\end{lemma}

Hence, an upper bound of $\bar{\Phi}$ can be obtained by minimizing $C$ over $V$ under the constraint $H(\bfx)\le C$,
which can be formulated as a SOS optimization problem in the form:
\bea
&\displaystyle\min_{V}~C \label{SOS_1}\\
[1ex]
&\mbox{s.t.} ~~-\left( \bff(\bfx) \cdot\nabla_{\bfx} V(\bfx)+\Phi(\bfx)-C\right) \mbox{~is~SOS},
\label{SOS}
\eea
which is a special case of (\ref{linear}). A better bound might be obtained by removing the requirement for $V(\bfx)$ to be a  polynomial and replacing (\ref{SOS}) with the requirement of non-negativeness.
However, the resulting problem would be too difficult, since the classical algebraic-geometry problem
of verifying positive-definiteness of a general multi-variate polynomial is NP-hard \cite{An:02,An:05}.

Notice that while  $V$ is similar to a Lyapunov function in a stability analysis, it is not required to be positive-definite. Notice also that a lower bound of any long-time average cost of the system (\ref{sys1}) can be analyzed in a similar way.

\begin{remark}\label{RemarkOnBoundedness}
For many systems the boundedness of system state immediately follows from energy consideration. In general, if the system state is bounded, this can often be proven using the SOS approach.
It suffices to check whether there exists a large but bounded global attractor, denoted by $\mathcal{D}_1.$
As an example, let $\mathcal{D}_1=\{\bfx~|~0.5\bfx^T\bfx\le \beta\}$, where the constant $\beta$ is sufficiently large. Then, the global attraction property of system in $\mathcal{D}_1$ may be expressed as
\bea
\bfx^T\dot{\bfx}=\bfx^T\bff(\bfx)\le -(0.5\bfx^T\bfx-\beta).
\label{SOS1}
\eea
Introducing a tunable polynomial $S(\bfx)$ satisfying $S(\bfx)\ge 0 ~\forall \bfx\in{\mathbb R}^n$, by Lemma~1, (\ref{SOS1}) can be relaxed to
\bea
\left\{
\begin{array}{c}
-\left(\bfx^T\bff(\bfx)-S(\bfx)(0.5\bfx^T\bfx-\beta)\right)\mbox{~is ~SOS},  \\
S(\bfx)\mbox{~is ~SOS}.
\end{array}
\right.
\label{SOS2}
\eea
Minimization of upper bound of long-time average cost for systems that have unbounded global attractor is usually meaningless, since the cost itself could be infinitely large.
\end{remark}


\section{Problem Formulation}\label{seq:Problem Formulation}

Consider a polynomial system with single input
\bea
\dot{\bfx}=\bff(\bfx)+\bfg(\bfx)\bfu
\label{sys}
\eea
where $\bff(\bfx): {\mathbb R}^n\rightarrow {\mathbb R}^n$ and $\bfg(\bfx): {\mathbb R}^n\rightarrow {\mathbb R}^{n\times m}$ are
polynomial functions of system state $\bfx$. The approach of this paper can easily be extended to multiple input systems.
The control $\bfu\in {\mathbb R}^m$, which is assumed to be a polynomial vector of the system state $\bfx$ with maximum degree $d_{\bfu}$,
is designed to minimize the upper bound of an average cost of the form:
\bea
\bar{\Phi}=\lim_{T\rightarrow \infty}\frac{1}{T}\int_0^T\Phi({\bfx}(t),\bfu(t))\,dt,
\label{cost}
\eea
where ${\bfx}$ is the closed-loop solution of the system (\ref{sys}) with the control $\bfu$.
The continuous function $\Phi$ is a given non-negative polynomial cost in ${\bfx}$ and $\bfu$.

Similarly to (\ref{SOS_1})-(\ref{SOS}), we consider the following optimization problem:
\bea
&\displaystyle \min_{\bfu, V} C \label{objective}\\
&{s.t.} -\left((\bff(\bfx)+\bfg(\bfx)\bfu)\cdot\nabla_{\bfx} V+\Phi(\bfx, \bfu)-C\right) \mbox{is SOS}.\quad
\label{optimization}
\eea
When it cannot be guaranteed that the closed-loop system state is bounded, SOS constraints (\ref{SOS2}) must be added to (\ref{optimization}) to make
our analysis rigorous.

Under the framework of SOS optimization, the main problem in solving (\ref{objective})-(\ref{optimization}) is due to the non-convexity of (\ref{optimization}) caused by the control input $u$ and the decision function $V,$ both of which are tunable, entering~(\ref{optimization}) nonlinearly.
Iterative methods \cite{Zh:07,Zh:09,Ng:11} may help to overcome this issue indirectly in the following way: first fix one subset of bilinear decision
variables and solve the resulting linear inequalities in the other decision variables; in the next step, the other bilinear decision variables are fixed and the procedure is repeated. For the particular long-time average cost control problem (\ref{objective})-(\ref{optimization}),
the non-convexity will be resolved in the following by considering a type of so-called small-feedback controller.
In such a new way, iterative updating of decision variables is exempted, and replaced by solving a sequence of SOS optimization problems.


\section{Bound optimization of long-time average cost for controlled polynomial systems}\label{seq:Bound}

In this section a small-feedback controller is designed to reduce the upper bound of the long-time average cost (\ref{cost}) for the controlled polynomial system (\ref{sys}).
It is reasonable to hope that a controller reducing the upper bound for the time-averaged cost will also reduce the time-averaged cost itself \cite{Ph:14}.

\subsection{Basic formalism of the controller design}

We will look for a controller in the form
\bea
\bfu(\bfx,\epsilon)=\sum_{i=1}^{\infty}\epsilon^i \bfu_i(\bfx),
\label{cc}
\eea
where $\epsilon>0$ is a parameter, and $\bfu_i(\bfx), i= 1,2,\cdots$ are polynomial vector functions of system state $\bfx.$
In other words, we seek a family of controllers parameterised by $\epsilon$ in the form of a Taylor series in $\epsilon$. Notice that the expansion starts at the first-order term, so that $\epsilon=0$ gives the uncontrolled system.
To resolve the non-convexity problem of SOS optimization, we expand $V$  and $C$ in $\epsilon$:
\bea
      V(\bfx,\epsilon)&=&\sum_{i=0}^{\infty} \epsilon^i V_i(\bfx), \label{LF0} \\
      C(\epsilon)&=&\sum_{i=0}^{\infty}\epsilon^i C_i,
\label{LF1}
\eea
where $V_i$ and $C_i$ are the Taylor series coefficients for the tunable function and the bound, respectively, in the $i$th-order term of $\epsilon$.
Define
\bea
F(V,u,C)\eqdef (\bff(\bfx)+\bfg(\bfx)\bfu)\cdot \nabla_{\bfx} V+\Phi(\bfx,\bfu)-C.
\label{cc1}
\eea
Substituting (\ref{cc}), (\ref{LF0}), and (\ref{LF1}) into (\ref{cc1}), we have
\beas
F(V,u,C)&=& \left(\bff+\bfg\sum_{i=1}^{\infty}\epsilon^i\bfu_i\right)\cdot \sum_{i=0}^{\infty}\epsilon^i\nabla_{\bfx} V_i+\Phi\left(\bfx,\sum_{i=1}^{\infty}\epsilon^i\bfu_i\right)
-\sum_{i=0}^{\infty}\epsilon^i C_i.
\eeas
Noticing
\beas
\Phi\left(\bfx,\sum_{i=1}^{\infty}\epsilon^i\bfu_i\right)
=
\sum_{i=0}^{\infty}\epsilon^i\left(\sum_{k=0}^i\frac{1}{k!}\frac{\partial^k\Phi}{\partial \bfu^k}(\bfx,0)\frac{1}{i!}\frac{\partial^i\left(\bfu^k\right)}{\partial \epsilon^i}(\bfx,0)\right),
\eeas
it follows that
\bea
F(V,u,C)=\sum_{i=0}^{\infty}\epsilon^i F_i(V_0,\cdots,V_i, \bfu_{1},\cdots,\bfu_i, C_i),
\label{inq2}
\eea
where
\bea
F_i=\bff\cdot \nabla_{\bfx}V_i+\sum_{j+l=i}\bfg\bfu_j\cdot\nabla_{\bfx}V_l
+\sum_{k=0}^i\frac{1}{k!}\frac{\partial^k\Phi}{\partial \bfu^k}(\bfx,0)\frac{1}{i!}\frac{\partial^i\left(\bfu^k\right)}{\partial \epsilon^i}(\bfx,0)-C_i.
\label{new1}
\eea
In (\ref{new1}), $({\partial^k\Phi}/{\partial \bfu^k})(\bfx,0)$ denotes the $k$th partial derivative of $\Phi$ with respect to $\bfu$ at $\bfu=0$, and
$({\partial^i(\bfu^k)}/{\partial \epsilon^i})(\bfx,0)$ denotes the $i$th partial derivative of $\bfu^k(\bfx,\epsilon)\eqdef[u_1^k(\bfx,\epsilon),\cdots, u_m^k(\bfx,\epsilon)]^T$ with respect to $\epsilon$
at $\epsilon=0$.

Expression (\ref{inq2}) becomes more clear when a specific cost function $\Phi$ is considered.
For instance, let $\Phi=\Phi_0(\bfx)+\bfu^T\bfu$.
Then,
\beas
F(V,\bfu,C)
=F_0(V_0,C_0)+\epsilon F_1(V_0,V_1,\bfu_1,C_1)+\epsilon^2F_2(V_0,V_1,V_2, \bfu_1,\bfu_2,C_2)+O(\epsilon^3),
\eeas
where
\beas
F_0&=&\bff\cdot \nabla_{\bfx} V_0+\Phi_0-C_0, \\
[1ex]
F_1&=& \bff\cdot \nabla_{\bfx} V_1+\bfg \bfu_1\cdot \nabla_{\bfx} V_0-C_1,\\
[1ex]
F_2&=&\bff\cdot \nabla_{\bfx} V_2+\bfg \bfu_1\cdot \nabla_{\bfx} V_1+\bfg \bfu_2\cdot \nabla_{\bfx} V_0+\bfu_1^T\bfu_1-C_2,
\eeas
and $O(\epsilon^3)$ denotes all the terms with order of $\epsilon$ being equal or greater than 3.

It is clear that $F(V,\bfu,C)\le 0$ holds if $F_i\le 0, i=0,1,2,\cdots$, simultaneously, and the series (\ref{cc})-(\ref{LF1}) converge.
Notice that $F_i$ includes tunable functions $V_j, j\le i$, and $\bfu_k, k\le i-1$. For any non-negative integers $i_1, i_2$ satisfying $i_1<i_2$, the tunable variables in $F_{i_1}$ are always a subset of the tunable variables in $F_{i_2}$. Hence (\ref{objective})-(\ref{optimization}) can be solved as a sequence of convex optimization problems. When the inequality constraints $F_i\le 0$ are relaxed to SOS conditions, our idea can be summarized as follows.

\noindent\rule{16.5cm}{0.1pt}
{\bf\it The sequential steps to solve (\ref{objective})-(\ref{optimization}): {\bf A-I}} \\
\noindent\rule{16.5cm}{0.1pt}

\begin{itemize}

\item[(s0)] First minimize $C_0$ over $V_0$ under the constraint $F_0(V_0,C_0)\le 0$, or more conservatively,
\beas
O_0: ~~\min_{V_0} C_0, {~~s.t.~~} -F_0(V_0,C_0)\mbox{~~is SOS}.
\eeas
Denote the optimal $C_0$ by $C_{0,SOS}$ and the associated $V_0$ by $V_{0,SOS}$.

\item[(s1)] Now, let $V_0=V_{0,SOS}$ in $F_1$, and then minimize $C_1$ over $V_1$ and $\bfu_1$ under the constraint $F_1(V_{0,SOS},V_1,\bfu_1,C_1)\le 0$, or under the framework of SOS optimization,
\beas
O_1: ~~\min_{V_1,\bfu_1} C_1, {~s.t.~} -F_1(V_{0,SOS},V_1,\bfu_1,C_1)\mbox{~is~ SOS}.
\eeas
Using the generalized $\mathcal{S}$-procedure given in Lemma~1 and the fact that
\bea
-F_0(V_{0,SOS},C_{0,SOS})\ge 0,
\label{condition1}
\eea
$O_1$ can be revised by incorporating one more tunable function $S_0(\bfx)$:
\beas
O_1': \begin{array}{c}
\min_{V_1,\bfu_1,S_0} C_1, \\
[1ex]
{~~s.t.~~} \left\{
\begin{array}{c}
-F_1(V_{0,SOS},V_1,\bfu_1,C_1)
+S_0(\bfx)F_0(V_{0,SOS},C_{0,SOS}) \mbox{~~is SOS}, \\
[1ex]
S_0(\bfx)\mbox{~~is~ SOS}.
\end{array}
\right.
\end{array}
\eeas
Denote the optimal $C_1$ by $C_{1,SOS}$ and the associated $V_1$ and $\bfu_1$ by $V_{1,SOS}$ and $\bfu_{1,SOS}$, respectively.

\item[(s2)] Further let $V_0=V_{0,SOS}$, $V_1=V_{1,SOS}$, and $\bfu_1=\bfu_{1,SOS}$ in $F_2$, and then minimize $C_2$ over $V_2$ and $\bfu_2$ under the constraint $F_2(V_{0,SOS},V_{1,SOS},V_2, \bfu_{1,SOS}, \bfu_2,C_2)\le 0$. In a more tractable way, consider
\beas
O_2:
\begin{array}{c}
 \displaystyle \min_{V_2,~ \bfu_2} C_2,~ {s.t.}\\
-F_2(V_{0,SOS},V_{1,SOS},V_2, \bfu_{1,SOS}, \bfu_2,C_2)\mbox{~is SOS}.
\end{array}
\eeas
Similarly as in (s1), noticing (\ref{condition1}) and
\beas
-F_1(V_{0,SOS},V_{1,SOS},\bfu_{1,SOS},C_{1,SOS})\ge 0,
\eeas
the SDP problem $O_2$ can be revised by the generalized $\mathcal{S}$-procedure to the following form:
\beas
O_2': \begin{array}{c}
 \min_{V_2,\bfu_2,S_0,S_1} C_2,  ~~~{s.t.} \\
 [1ex]
\left\{
\begin{array}{c}
 -F_2(V_{0,SOS},V_{1,SOS},V_2, \bfu_{1,SOS}, \bfu_2,C_2)
 +S_0(\bfx)F_0(V_{0,SOS},C_{0,SOS}) \\
~ +S_1(\bfx)F_1(V_{0,SOS},V_{1,SOS}, \bfu_{1,SOS},C_{1,SOS})\mbox{~is SOS},\\
 [1ex]
 S_0(\bfx)\mbox{~~is SOS}, \\
 [1ex]
 S_1(\bfx)\mbox{~~is SOS}.
\end{array}
\right.
\end{array}
\eeas
Denote the optimal $C_2$ by $C_{2,SOS}$ and the associated $V_2$ and $\bfu_2$ by $V_{2,SOS}$ and $\bfu_{2,SOS}$, respectively.

Notice that $S_0(\bfx)$ here might differ from the tunable function $S_0(\bfx)$ in $O_1'$. Throughout this paper we will use the same notations for the tunable functions like $S_0$ and $S_1$ in various instances of the  $\mathcal{S}$-procedure, to keep the notation simple.

\item[(s3)] The SOS-based controller design procedure is continued for higher-order terms.

\end{itemize}
\noindent\rule{16.5cm}{0.1pt}

Now, define three series
\bea
C_{SOS}=\sum_{i=0}^{\infty}\epsilon^iC_{i,SOS}, ~~\bfu_{SOS}=\sum_{i=1}^{\infty}\epsilon^i \bfu_{i,SOS}, ~~V_{SOS}=\sum_{i=0}^{\infty}\epsilon^iV_{i,SOS}.
\label{series}
\eea
When all of them converge, the following statement will be true.

\begin{theorem}
By applying the state-feedback controller $\bfu=\bfu_{SOS}$ for the system (\ref{sys}), if the trajectories of the closed-loop system are bounded
\footnote{In the context of long-time average cost controller design and analysis, it is actually enough to assume the boundedness of the global attractor of the system to ensure the existence of $C_{SOS}$.
}, then $C_{SOS}$ is an upper bound of the long-time average cost $\bar{\Phi}$.
\end{theorem}

{\it Proof}.  Using the algorithm {\bf A-I}, we obtain
\beas
F_i(V_{0,SOS},\cdots,V_{i,SOS}, \bfu_{1,SOS},\cdots, \bfu_{i,SOS},C_{i,SOS})\le 0, \forall ~i.
\label{new2}
\eeas
Then, it follows that
\beas
\sum_{i=0}^{\infty}F_i(V_{0,SOS},\cdots,V_{i,SOS}, \bfu_{1,SOS},\cdots, \bfu_{i,SOS},C_{i,SOS})
=F(V_{SOS}, \bfu_{SOS}, C_{SOS})\le 0,
\label{new3}
\eeas
where $C_{SOS}, \bfu_{SOS}, V_{SOS}$ are given in (\ref{series}).
By virtue of a same analysis as in proving Lemma~2 (see \cite{Ph:14}), we can conclude that $\bar{\Phi}\le C_{SOS}$.
\rule{0.09in}{0.09in}

\begin{remark}
After specifying the structure of controller to be of the form (\ref{cc}), the non-convexity in solving the optimization problem (\ref{objective})-(\ref{optimization}) has been avoided by solving the linear SDPs $O_0, O_1', O_2', \cdots$ in sequence.
During the process, all the involved decision variables are optimized sequentially, but not iteratively as in other methods \cite{Zh:07,Zh:09,Ng:11}.
\end{remark}


\begin{remark}
The smallness of $\epsilon$ can be used to relax $O_1', O_2', \cdots$ further.
For instance, in $O_1'$, in order to prove
$
F_0+\epsilon F_1\le 0,
$
we prove $F_1(V_{0,SOS},V_1,\bfu_1,C_1)\le 0$ with the aid of the known constraint $F_0(V_{0,SOS},C_{0,SOS})\le 0$,
thus not using that $\epsilon$ is small.
In fact, when $\epsilon$ is small, for $F_0+\epsilon F_1$ to be negative $F_1$ has to be negative only for those $\bfx$ where $F_0(\bfx)$ is small,
and not for all $\bfx$ as required in $O_1'$.
Meanwhile, checking the convergence of the series (\ref{series}) would be challenging or even impractical.
These points will be addressed in what follows.
\end{remark}

\subsection{Design of small-feedback controller}

Next, the sequential design method {\bf A-I} is revised to utilize that $\epsilon\ll 1$.

\noindent\rule{16.5cm}{0.1pt}
{\bf\it The revised sequential steps to solve (\ref{objective})-(\ref{optimization}): {\bf A-II}}

\noindent\rule{16.5cm}{0.1pt}

\begin{itemize}

\item[(s0)] Same as in {\bf A-I}, first solve the SOS optimization problem $O_0$.
Denote the optimal $C_0$ by $C_{0,SOS}$ and the associated $V_0$ by $V_{0,SOS}$.

\item[(s1)] Let $V_0=V_{0,SOS}$ in $F_1$, and then consider the following SDP problem:
\beas
&O_1'': \begin{array}{c}
\displaystyle \min_{V_1,\bfu_1,S_0} C_1, \\
[1ex]
{~~s.t.~~}
\begin{array}{c}
-F_1(V_{0,SOS},V_1,\bfu_1,C_1)
+S_0(\bfx)F_0(V_{0,SOS},C_{0,SOS}) \mbox{~is SOS},
\end{array}
\end{array}
\eeas
where $S_0$ is any tunable polynomial function of $\bfx$ of fixed degree.
Denote the optimal $C_1$ by $C_{1,SOS}$ and the associated $V_1$ and $\bfu_1$ by $V_{1,SOS}$ and $\bfu_{1,SOS}$, respectively.
Unlike $O_1'$, here the non-negativity requirement of $S_0$ is not imposed.
This can be understood as that the non-negativity constraint  is imposed only for $\bfx$ such that $F_0(V_{0,SOS},C_{0,SOS})=0$.


\item[(s2)] Further let $V_0=V_{0,SOS}$, $V_1=V_{1,SOS}$, and $\bfu_1=\bfu_{1,SOS}$ in $F_2$, and then consider
\beas
O_2'': \begin{array}{c}
\displaystyle  \min_{V_2,\bfu_2,S_0,S_1} C_2, ~~~{s.t.} \\
 [1ex]
\left\{
\begin{array}{c}
-F_2(V_{0,SOS},V_{1,SOS},V_2, \bfu_{1,SOS}, \bfu_2,C_2)
 +S_0(\bfx)F_0(V_{0,SOS},C_{0,SOS}) \\
~+S_1(\bfx) F_1(V_{0,SOS},V_{1,SOS},\bfu_{1,SOS},C_{1,SOS})
\mbox{~~is~ SOS},
\end{array}
\right.
\end{array}
\eeas
where $S_0$ and $S_1$ are any tunable polynomial functions of fixed degrees.
$S_0$ here does not need to be the same as in $O_1''$.
Denote the optimal $C_2$ by $C_{2,SOS}$ and the associated $V_2$ and $\bfu_2$ by $V_{2,SOS}$ and $\bfu_{2,SOS}$, respectively.
Similarly as in $O_1''$, here the non-negativity constraint is in effect imposed only where $F_0(V_{0,SOS},C_{0,SOS})=F_1(V_{0,SOS},V_{1,SOS}, \bfu_{1,SOS},C_{1,SOS})=0$.

\item[(s3)] The revised SOS-based controller design procedure is continued for higher-order terms.
\end{itemize}
\noindent\rule{16.5cm}{0.1pt}

Since the constraints of $S_i$ being SOS imposed in {\bf A-I} are removed in {\bf A-II}, the coefficients $C_{i,SOS}$ obtained in {\bf A-II} can be smaller than the coefficients $C_{i,SOS}$ obtained in {\bf A-I}. This advantage comes at a price: even if all the relevant series converge for a particular value of $\epsilon$, the procedure {\bf A-II} does not guarantee that the value $C_{SOS}$ given in (\ref{series}) is an upper bound for the time-averaged cost of the closed-loop system with the controller $\bfu_{SOS}$. We have now to consider (\ref{series}) as asymptotic expansions rather than Taylor series. Accordingly, we have to truncate the series and hope that the resulting controller will work for (sufficiently) small $\epsilon$
\footnote{It is worthy of noticing that the series truncation here does not mean that our controller design and analysis are conducted in a non-rigorous way.
The truncated controller would be effective if it leads to a better (lower) bound of the long-time average cost.}. It is possible to prove that this is, indeed, the case.

For illustration, the first-order truncation is considered only.


\begin{theorem}
Consider the first-order small-feedback controller for the system (\ref{sys}),
\bea
\bfu_{SOS}=\epsilon \bfu_{1,SOS}
\label{new3}
\eea
where $\epsilon>0$ is sufficiently small. Assume that the trajectories of the closed-loop system are bounded, and that $C_{1,SOS}<0$.
Then, $C_{\kappa,SOS}\eqdef C_{0,SOS}+\epsilon\kappa C_{1,SOS}, \kappa\in(0,1)$ is an upper bound of the long-time average cost $\bar{\Phi}$.
Clearly, $C_{\kappa,SOS}<C_{0,SOS}$.
\end{theorem}

{\it Proof}. Let $V_{SOS}=V_{0,SOS}+\epsilon V_{1,SOS}$.
By substituting $V=V_{SOS}, C=C_{\kappa, SOS}, \bfu=\bfu_{SOS}$ in the constraint function $F(V,\bfu,C)$ that is defined in (\ref{cc1}), the remaining task
is to seek small $\epsilon>0$ such that
\bea
F(V_{SOS},\bfu_{SOS},C_{\kappa, SOS})
\le
0.
\label{cc1_1}
\eea
Notice that
\bea
F(V_{SOS},\bfu_{SOS},C_{\kappa,SOS})
=
F_0(V_{0,SOS},C_{0,SOS})
+\epsilon F_1(V_{0,SOS},V_{1,SOS},\bfu_{1,SOS},C_{1,SOS})
+\epsilon (1-\kappa)C_{1,SOS}+\epsilon^2 w(\bfx,\epsilon),
\label{boundx}
\eea
where
\beas
w(\bfx,\epsilon)=\bfg\bfu_1\cdot  \nabla_{\bfx}V_{1,SOS}
+\frac{1}{\epsilon^2}\left(\Phi(\bfx,\epsilon \bfu_{1,SOS})-\Phi(\bfx,0)-\epsilon\frac{\partial \Phi}{\partial \bfu}(\bfx,0)\bfu_{1,SOS}\right),
\eeas
and $F_0, F_1$, being polynomial in $\bfx$, possess all the continuity properties implied by the proof.
Let $\mathcal{D}\in{\mathbb R}^n$ be the phase domain that interests us, where the closed-loop trajectories are all bounded.
Then,
\bea
F_{1,max}\eqdef \max_{\bfx\in\mathcal{D}}F_1(V_{0,SOS},V_{1,SOS},\bfu_{1,SOS},C_{1,SOS})<\infty,
\label{bound1}
\eea
and $w(\bfx,\epsilon)$ is bounded for any $\bfx\in\mathcal{D}$ and any finite $\epsilon$ (the latter following from the standard mean-value-theorem-based formula for the Lagrange remainder). By (\ref{boundx}) and (\ref{bound1}),
\bea
F(V_{SOS},\bfu_{SOS},C_{\kappa,SOS})\le
F_0(V_{0,SOS},C_{0,SOS})+\epsilon F_{1,max}+\epsilon (1-\kappa)C_{1,SOS}+O(\epsilon^2).
\label{bound1_ex}
\eea
Meanwhile, consider the two inequality constraints obtained by solving $O_0$ and $O_1''$:
\bea
\left\{
\begin{array}{c}
 F_0(V_{0,SOS},C_{0,SOS})\le 0, \\
[1ex]
 F_1(V_{0,SOS},V_{1,SOS}, \bfu_{1,SOS},C_{1,SOS})\le 0
~~\forall \bfx ~~\mbox{such that}~ F_0(V_{0,SOS},C_{0,SOS})= 0.
\end{array}
\right.
\label{cons}
\eea

Define $\mathcal{D}_{\delta}\eqdef \left\{\bfx\in \mathcal{D} ~|~ \delta\le F_0(V_{0,SOS},C_{0,SOS})\le 0\right\}$ for a given constant $\delta\le 0$.
Clearly, $\mathcal{D}_{\delta}\rightarrow \mathcal{D}_{0}$ as $\delta\rightarrow 0$.
Further define
\bea
F_{1,\delta}(\delta)\eqdef \max_{\bfx\in \mathcal{D}_{\delta}} F_1(V_{0,SOS},V_{1,SOS},\bfu_{1,SOS},C_{1,SOS}).
\label{bound2}
\eea
By the second constraint in (\ref{cons}), $\lim_{\delta\rightarrow 0}F_{1,\delta}(\delta)\le 0$.
Therefore, by continuity and the fact $C_{1,SOS}<0$, for any $0<\kappa<1$ there exists a constant $\delta_{\kappa}<0$ such that
\bea
F_1(V_{0,SOS},V_{1,SOS}, \bfu_{1,SOS},C_{1,SOS})\le F_{1,\delta_{\kappa}}<-\frac{1}{2}(1-\kappa)C_{1,SOS}, ~~\forall \bfx\in \mathcal{D}_{\delta_{\kappa}}.
\label{bound3}
\eea
In consequence, (\ref{boundx}), the first constraint in (\ref{cons}), and (\ref{bound3}) render to
\bea
F(V_{SOS},\bfu_{SOS},C_{\kappa,SOS}) &\le&
F_0(V_{0,SOS},C_{0,SOS})+\epsilon F_{1,\delta_{\kappa}}+\epsilon (1-\kappa)C_{1,SOS}+O(\epsilon^2) \nonumber \\
&\le& \frac{\epsilon}{2} (1-\kappa)C_{1,SOS}+O(\epsilon^2)\le 0, ~~\forall \bfx\in \mathcal{D}_{\delta_{\kappa}},
\label{bound4}
\eea
for sufficiently small $\epsilon$.

Next, we prove (\ref{cc1_1}) for any $\bfx\in \mathcal{D}\setminus\mathcal{D}_{\delta_{\kappa}}$. By the definition of the set $\mathcal{D}_{\delta_{\kappa}}$, we have
\bea
F_0(V_{0,SOS},C_{0,SOS})<\delta_{\kappa}<0,  ~~\forall \bfx\in \mathcal{D}\setminus\mathcal{D}_{\delta_{\kappa}}.
\label{bound5}
\eea
Then, (\ref{bound1_ex}) and (\ref{bound5}) yield
\bea
F(V_{SOS},\bfu_{SOS},C_{\kappa,SOS})\le
\delta_{\kappa}+\epsilon F_{1,max}+\epsilon (1-\kappa)C_{1,SOS}+O(\epsilon^2)\le \delta_{\kappa}+O(\epsilon)\le 0, ~~\forall \bfx\in \mathcal{D}\setminus\mathcal{D}_{\delta_{\kappa}},
\label{bound1_ex1}
\eea
if $\epsilon$ is sufficiently small.

(\ref{bound4}) and (\ref{bound1_ex1}) imply that (\ref{cc1_1}) holds $\forall \bfx\in \mathcal{D}$. The proof is complete.
\rule{0.09in}{0.09in}


In practice, once the form of the controller has been specified in (\ref{new3}), the upper bound $C$ and the corresponding $V$ actually can be obtained by solving the following optimization problem directly:
\beas
O_{\epsilon}: \begin{array}{c}
\displaystyle  \min_{V,~\epsilon} C, ~~~ \\
 [1ex]
{s.t.}~~-F(V,\epsilon \bfu_{1,SOS},C)\mbox{~is~ SOS}.
\end{array}
\eeas
This problem can be further relaxed by incorporating the known constraints (\ref{cons}).
In $O_{\epsilon}$, if $\epsilon$ is set as one of the tunable variables, the SOS optimization problem will become non-convex again, thus causing additional trouble in solving it. Alternatively, one can fix $\epsilon$ here, and investigate its effect on the upper bound of $\bar{\Phi}$ by trial and error. We will follow this route in Section~\ref{seq:Example}.

\section{Illustrative example}\label{seq:Example}

As an illustrative example we consider a system proposed in \cite{Ki:05} as a model for studying control of oscillatory vortex shedding behind a cylinder. The actuation was assumed to be achieved by a volume force applied in a compact support region downstream of the cylinder. The Karhunen-Lo\`{e}ve (KL) decomposition \cite{No:03} was used and the first two KL modes and an additional shift mode were selected. For the Reynolds number equal to 100 the resulting low-order Galerkin model of the cylinder flow with control was given as follows
\bea
\left[
\begin{array}{c}
\dot{a}_1 \\
\dot{a}_2 \\
\dot{a}_3
\end{array}
\right]&=&
\left[
\begin{array}{ccc}
\sigma_r & -\omega-\gamma a_3 & -\beta a_1 \\
\omega+\gamma a_3 & \sigma_r & -\beta a_2\\
\alpha a_1 & \alpha a_2 & -\sigma_3
\end{array}
\right]
\left[
\begin{array}{c}
{a}_1 \\
{a}_2 \\
{a}_3
\end{array}
\right]
+
\left[
\begin{array}{c}
g_1 \\
g_2 \\
0
\end{array}
\right]u,
\label{cf}
\eea
where $\sigma_r=0.05439, \sigma_3=0.05347, \alpha=0.02095, \beta=0.02116,$
$\gamma=-0.03504, \omega=0.9232, g_1=-0.15402$, and $g_2=0.046387$.
More details on deriving the reduced-order model (\ref{cf}) are given in~\cite{Ro:14}.

The system (\ref{cf}) possesses a unique equilibrium when $u=0$, which is at the origin.
Let $\Phi=1/2\bfa^T\bfa+u^2$, where $\bfa=[a_1 ~a_2 ~a_3]^T$.
The proposed algorithms {\bf A-I} and {\bf A-II} were applied to (\ref{cf}), with
the system state assumed to be available. In experiment, it could be estimated by designing a state observer with some sensed output measurement at a typical position~\cite{Ro:14}.

\subsection{Performance of algorithm {\bf A-I}}


The SDP problem $O_0$ is solved first. It corresponds to the uncontrolled sysytem.
The minimal upper bound we could achieve was $C_{0,SOS}=6.59.$ It was obtained with
\beas
V_{0,SOS}=-96.63 a_3+14.01a_1^2+14.01a_2^2+14.15 a_3^2.
\eeas
Increasing the degree of $V_0$ cannot give a better bound because there exists a stable limit cycle in the phase space of (\ref{cf}), on which $a_1^2+a_2^2=6.560,$ and $ a_3=2.570$.
Since $\bar{\Phi}=1/2\bfa^T\bfa=6.584$ on the limit cycle, the minimal upper bound achieved by SOS optimization is tight in the sense
that the difference between $C_{0,SOS}$ and $\bar{\Phi}$ is less than the prescribed precision for $C$, $0.01$.

Solving the SDP problem $O_1$, where $V_1$ and $u_1$ are tunable functions, gave $C_{1,SOS}=0$.
Solving $O_1',$ with $V_1, u_1, S_0$ being tuning functions, gave the same result: $C_{1,SOS}=0$.
In both cases, increasing the degrees of the tuning functions did not reduce the upper bound.
The consequent SOS optimization problems, $O_i',$ with $i=2,3$ also gave $C_{i,sos}=0, i=2,3$.
Therefore, by (\ref{series}),
\beas
C_{SOS}=C_{0,SOS}+\epsilon C_{1,SOS}+\epsilon^2 C_{2,SOS}+O(\epsilon^3)
       \approx C_{0,SOS}=6.59,
\eeas
implying that {\bf A-I} does not generate a control ensuring a better upper bound of $\bar{\Phi}$ than the bound obtained in the uncontrolled case.

\subsection{Performance of algorithm {\bf A-II}}

Without any control, it has been obtained in {\bf A-I} that $C_{0,SOS}=6.59$.

We first solve $O_1''$. Given the vectors of monomials in $\bfx$ without repeated elements \cite{Ch:09}, $Z_i, i=1,\cdots, 3$, define $V_1=P_1^TZ_1$, $u_1=P_2^TZ_2$, and $S_0=P_3^TZ_3$, where the parametric vectors $P_i$, $i=1,\cdots, 3$ consist of tuning vector variables. The degrees of $V_1, u_1$ and $S_0$ are specified by the maximum degrees of monomials in $Z_i, i=1,\cdots, 3$, and denoted by $ d_{V_1}, d_{u_1}$, and $d_{S_0}$, respectively. Consider two subcases: $d_{V_1}=d_{u_1}=d_{S_0}=2$ and $d_{V_1}=d_{u_1}=d_{S_0}=4$.
For the former case, we have $C_{1,SOS}=-354$, induced by
\beas
u_{1,SOS,2}=45.37 a_1-28.47 a_2-142.76 a_2a_3+399.49 a_1a_3.
\eeas
For the latter case, we have $C_{1,SOS}=-1965$, induced by
\beas
u_{1,SOS,4}&=&233.08a_1-54.73a_2-67.61a_2a_3+218.56a_1a_3+717.28a_1^3+13.16a_1^2a_2 \\
&&+571.67a_1a_2^2-277.73a_2^3+466.61a_1a_3^2-141.41a_2a_3^2+230.53a_1^3a_3\\
&&+106.32a_1^2a_2a_3+220.19a_1a_2^2a_3-161.44a_2^3a_3+628.40a_1a_3^3-173.78a_2a_3^3.
\eeas

\begin{figure}
  \centerline{\includegraphics[trim=0mm 5mm 0mm 0mm,clip,width=10cm]{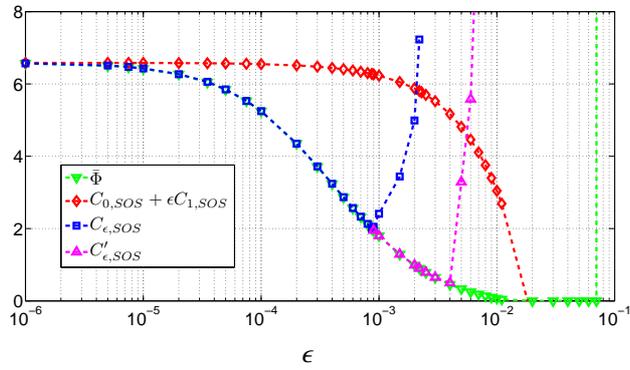}}
\centerline{$\epsilon$}
  \caption{The long-time average cost $\bar{\Phi}$ and its upper bounds for different $\epsilon$.
   $C_{0,SOS}, C_{1,SOS}, C_{\epsilon,SOS}, C_{\epsilon,SOS}'$ are obtained by solving $O_0, O_1'', O_{\epsilon}$, and $O_{\epsilon}$ with the relaxation (\ref{cons}), respectively.
}
\label{inf0}
\end{figure}

We then solve $O_{\epsilon}$ with a fixed $\epsilon$.
For simplicity we considered $u_{1,SOS}=u_{1,SOS,2}$ and $d_V\le 10$ only.
The upper-bound results for different $\epsilon$ are summarized in Fig. \ref{inf0}.
The long-time average cost $\bar{\Phi}$, which is obtained by direct numerical experiment, and the linear truncated bound $C_{0,SOS}+\epsilon C_{1,SOS}$ are also presented for comparison.
From Fig. \ref{inf0}, we can see the following.

Let $\epsilon_1= 1.267\times 10^{-2}$ and $\epsilon_2=7.416\times 10^{-2}$.
The small-feedback controller
\bea
u=\epsilon u_{1,SOS,2}
 \label{test}
\eea
 reduces $\bar{\Phi},$ and the reduction in $\bar{\Phi}$ increases monotonically with $\epsilon$ when $0<\epsilon<\epsilon_2$. In particular, $\bar{\Phi}=0$ for $\epsilon_1 \le \epsilon< \epsilon_2,$ that is in this range of $\epsilon$ the controller fully stabilizes the system. When $\epsilon\ge \epsilon_2$, the controller makes the long-time average cost worse than in the uncontrolled case. The effect of $\epsilon$ on $\bar{\Phi}$ can be seen more clearly by investigating the qualitative properties of the closed-loop system. A simple check gives that when $0\le \epsilon< \epsilon_1$, the closed-loop system has a unique unstable equilibrium at the origin and a stable limit cycle, thus yielding a non-zero but finite $\bar{\Phi}$;
when $\epsilon_1 \le \epsilon< \epsilon_2$, the limit cycle disappears and the unique equilibrium becomes globally stable, thus implying
the vanishness of $\bar{\Phi}$;
when $\epsilon\ge \epsilon_2$ but is close to $\epsilon_2$, besides the equilibrium at the origin, there exist four additional non-zero equilibria, and as a result $\bar{\Phi}$ becomes large immediately. For instance, at the bifurcation point $\epsilon=\epsilon_2$, the non-zero equilibria of the closed-loop system are $(\pm 0.6988, \pm 2.362, 2.377)$ and $(\pm 0.7000, \pm 2.364,  2.382)$, resulting in $\bar{\Phi}=171.55$.

Solving $O_{\epsilon}, 0<\epsilon\le 8.7\times 10^{-4}$ yields a tight upper bound $C_{\epsilon, SOS}$ for $\bar{\Phi}$.
However, the obtained upper bound becomes non-tight when $\epsilon> 8.7\times 10^{-4}$.
The conservativeness of $C_{\epsilon, SOS}$ can be fully overcome by considering additional relaxation constraint (\ref{cons})
for $ 8.7\times 10^{-4}<\epsilon\le 4\times 10^{-3}$, but only mitigated to certain extend for larger $\epsilon$.

The two-term expansion $C_{0,SOS}+\epsilon C_{1,SOS}$ is only a linear approximation of $C_{SOS}$ in (\ref{series}).
Thus, as an upper bound of $\bar{\Phi}$, it behaves well when $\epsilon$ is very small, but it becomes conservative when $\epsilon$ is further increased, and meaningless as $\epsilon>-C_{0,SOS}/C_{1,SOS}=0.0186$.

In summary, for small $\epsilon$, the proposed small-feedback controller
yields a better bound of the long-time average cost than in the uncontrolled case. Further, the controller indeed reduces the long-time average cost itself.


\begin{figure}
  \centerline {\includegraphics[trim=0mm 10mm 0mm 0mm,clip,width=10cm]{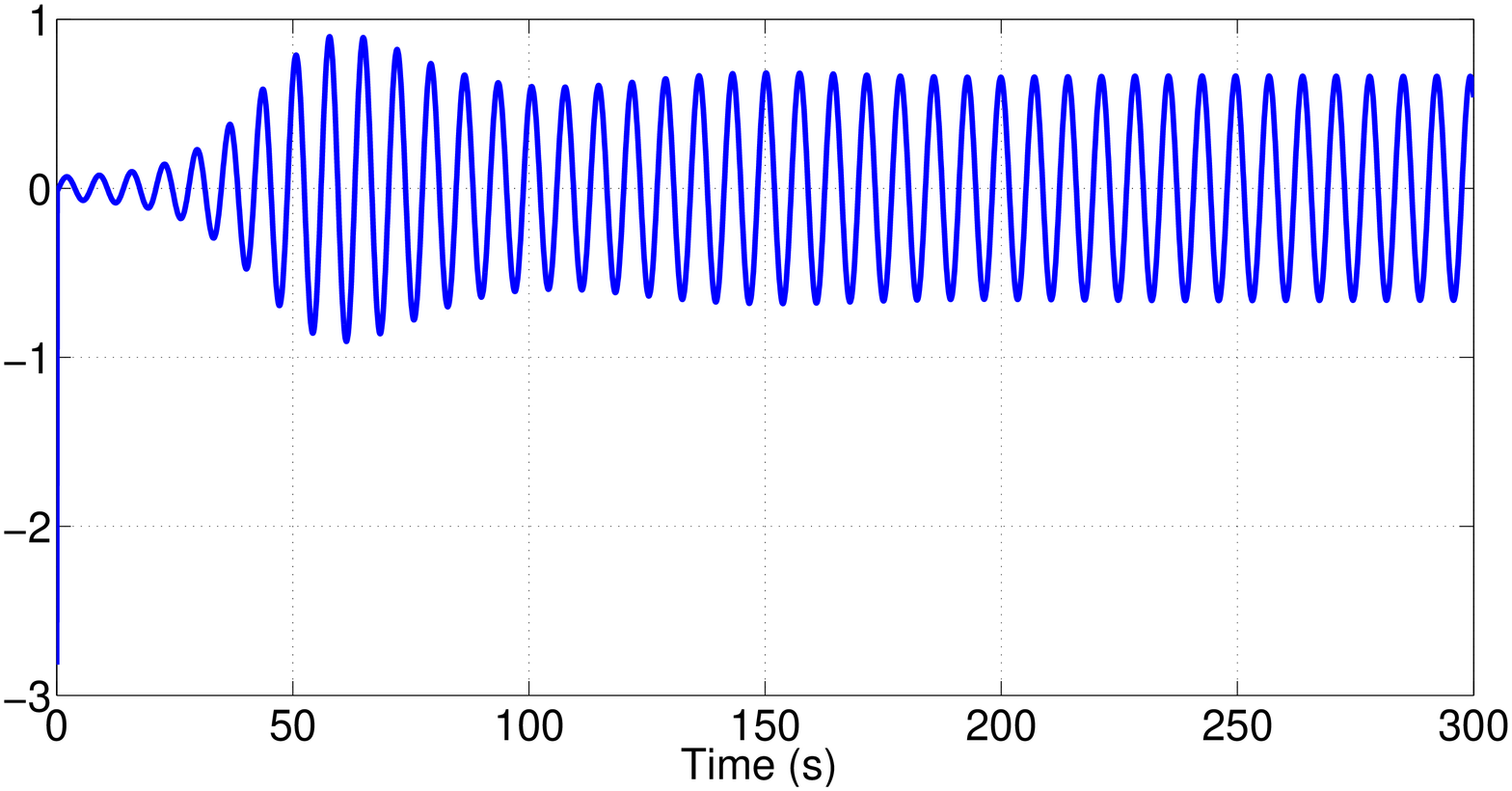}}

\centerline{$t$}
  \caption{Control input profile.
}
\label{inf2}
\end{figure}

\begin{figure}
  \centerline{\includegraphics[height=5cm,width=9cm]{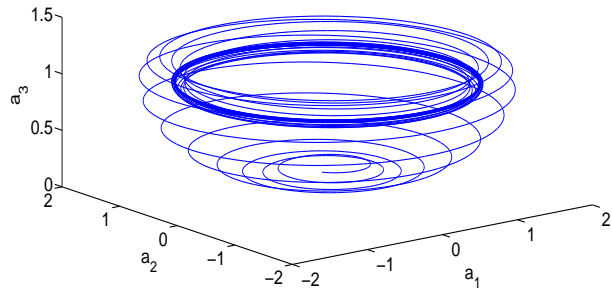}}
  \caption{Closed-loop trajectory starting at $\bfa=[-0.3~-0.3~0.3]^T$.
  Owing to the small-feedback control, the magnitude of the periodic oscillation has been reduced.
}
\label{inf3}
\end{figure}

Figs. \ref{inf2}-\ref{inf3} show more details of the control performance of the proposed controller (\ref{test}) with $\epsilon=8.7\times 10^{-4}$ and
the initial state $\bfa=[-0.3 ~-0.3 ~0.3]^T$.


\section{Conclusion}\label{seq:Conclusion}

Based on sum-of-squares decomposition of polynomials and semidefinite programming, a numerically tractable approach is presented for long-time average cost control of polynomial dynamical systems.
The obtained controller possesses a structure of small feedback, which is an asymptotic expansion in a small parameter, with all the coefficients being polynomials of the system state.
The derivation of the small-feedback controller is given in terms of the solvability conditions of state-dependent linear and bilinear inequalities.
The non-convexity in SOS optimization can be resolved by making full use of the smallness of the perturbation parameter while not using any iterative
algorithms.
The efficiency of the control scheme has been tested on a low-order model of  cylinder wake flow stabilization problem.
In the next research phase, we will consider SOS-based long-time average cost control under modelling uncertainties and in the presence of noise, as well as direct numerical simulations of small-feedback control for actual fluid flows.

The proof of concept of the idea of using the upper bound of long-time average cost control as the objective of the control design, and the method of overcoming the non-convexity of simultaneous optimization of the control law and the tunable function are the two main contributions of the present paper.

\end{document}